\def\goth{\mathfrak}

\documentstyle[amssymb,12pt]{amsart}

\newtheorem{theorem}{Theorem}
\newtheorem{application}{Application}

\newcommand{\res}{\upharpoonright}

\begin{document}

\baselineskip=20pt

  \begin{center}
     {\bf\large The Strength Of The Isomorphism Property}\footnote{{\em 
     Mathematics Subject Classification} Primary 03C50, 03H05. 
     Secondary 26E35, 28E05.}
  \end{center}

  \begin{center}
 Renling Jin\footnote{The first author would like to thank
 C. Ward Henson for some valuable suggestions} \& Saharon Shelah\footnote{The 
 research of the second author was partially
 supported by {\em Basic Research Fund, Israel Academy of Humanity
 and Sciences;} Pub. number 493.}
  \end{center}

  \bigskip

  \begin{quote}

    \centerline{Abstract}

    \small

In \S 1 of this paper, we characterize the isomorphism property 
of nonstandard universes in terms of
the realization of some second--order types in model theory. In \S 2, 
several applications are given. One of the applications
answers a question of D. Ross in [R] about infinite Loeb measure spaces.

\end{quote}

\section{Introduction}

We always use ${}^{*}V$ for a nonstandard universe. We refer to
[CK] or [SB] for the definition of nonstandard universes.

In the book [SB], there is an interesting
example (see [SB, Theorem 1.2.12.(e)]) for
illustrating the unusual behavior of infinite Loeb measure spaces.
The example of [SB] says that in a nonstandard universe called
a polyenlargement, the statement ($\dagger$) is true,
where the statement ($\dagger$) is the following: 

\begin{quote}

Every infinite Loeb measure space has a subset $S$ such that $S$
has infinite Loeb outer measure, but the intersection of $S$ with
any finite Loeb measure set has Loeb measure zero.

\end{quote} 
 
Under certain definition, the set $S$ is called measurable but has 
infinite outer measure and zero inner measure (see [SB]).
The diagonal argument for constructing $S$ in [SB] depends on the
construction of polyenlargements, say, an iterated
ultrapower (or ultralimit) construction.
During the preparation of the book [SB], K. D. Stroyan asked (see [R])
whether or not ($\dagger$) can be proved by some nice general properties 
of nonstandard universes without mentioning any particular construction.
The first natural candidate would be C. W. Henson's {\em isomorphism property}
[H1]. 

Let $\cal L$ be a first--order language. An $\cal L$--model $\goth A$
is called internally presented in $^{*}V$ if the base set $A$ and every
interpretation under $\goth A$ of a symbol in $\cal L$ are internal in $^{*}V$. 
(For any $\cal L$--model $\goth A$ and symbol $P$ in $\cal L$ we 
write $P^{\goth A}$ for the interpretation of $P$ under $\goth A$. We sometime
use ${\cal L}_{\goth A}$ for the language of $\goth A$.)
Let $\kappa$ be an infinite cardinal.
A nonstandard universe $^{*}V$ is said to satisfy the $\kappa$--isomorphism
property ($I\!P_{\kappa}$ for short) if

\begin{quote}

for any two internally presented $\cal L$--models $\goth A$ and $\goth B$
with $|{\cal L}|<\kappa$, ${\goth A}\equiv {\goth B}$ implies ${\goth A}\cong
{\goth B}$,

\end{quote}
where ``$\equiv$'' means {\em to be elementarily equivalent to} and ``$\cong$''
means {\em to be isomorphic to}. It is easy to see that $I\!P_{\kappa}$
implies $I\!P_{\kappa'}$ when $\kappa'\leq\kappa$.

Instead of using the isomorphism property, D. Ross in [R] proved that a
property called {\em the $\kappa$--special model axiom} for any infinite
cardinal $\kappa$, which is stronger than $I\!P_{\kappa}$, implies ($\dagger$). 
In [R], Ross showed also that 
${\kappa}$--special model axiom has many new consequences, which hadn't
been proved by $I\!P_{\kappa}$ then. In his paper, Ross asked which of those
results can or cannot be proved by $I\!P_{\kappa}$. The 
most important question among them is that if we can 
or cannot prove ($\dagger$) by $I\!P_{\kappa}$
for some infinite cardinal
$\kappa$. Basically, it was not known back then whether 
or not the $\kappa$--special model axiom is strictly stronger that 
$I\!P_{\kappa}$ (see [R]).

The first author then answered the most of Ross's questions in [J].
In that paper, Jin showed that $I\!P_{\kappa}$ for arbitrary large $\kappa$
does not imply some consequences of the $\aleph_0$--special model axiom. 
As a corollary $I\!P_{\kappa}$ is strictly weaker than the $\kappa$--special
model axiom.
He also showed that many of the consequences of the $\kappa$--special
model axiom in [R]
are also the consequences of $I\!P_{\kappa}$. Unfortunately, [J] didn't answer
Ross's question about ($\dagger$).

In the another direction, the authors of [JK] proved that ($\dagger$) is
true in some ultrapowers of the standard universe. Since 
we need iterated ultrapower construction to build the nonstandard
universes of the $\kappa$--special model axiom while we need only one--step
ultrapower construction to build
the nonstandard universes of $I\!P_{\kappa}$ (see [H2]), 
the result of [JK] seems to suggest that $I\!P_{\kappa}$ have the right
strength to prove ($\dagger$).

The main purpose of this paper is to solve Ross's question about ($\dagger$).
In \S 1, we characterize $I\!P_{\kappa}$ in terms of the realization of 
some second--order types. By applying Theorem 1 of \S 1, we show in \S 2 
that Ross's question about ($\dagger$) has a positive answer, {\em i.e.} 
($\dagger$) can be
proved by $I\!P_{\aleph_0}$. In \S 2, we reprove also three 
known results in [J] by using the same method in a uniform way. 
The new method simplifies significantly the original proofs in [J].

Notation for model theory in this paper will be consistent with [CK].

\section{Characterization of the isomorphism property}

We use always $\cal L$ for a first--order language. Let $X$ be an $n$--ary
predicate symbol which is not in $\cal L$. 
We call $\Gamma (X)$ an $n$--$\Delta^1_0 ({\cal L})$
type iff $\Gamma (X)$ is a consistent set of ${\cal L}\cup\{X\}$--sentences.
Let $\goth A$ be an $\cal L$--model with base set $A$ and let $\Gamma (X)$ be
an $n$--$\Delta^1_0 ({\cal L})$ type. We say that $\Gamma (X)$ is consistent
with $\goth A$ iff $\Gamma (X)\cup T\!h({\goth A})$ is consistent, where
$T\!h({\goth A})$ is the set of all $\cal L$--sentences which are true in 
$\goth A$. We say that $\goth A$ realizes $\Gamma (X)$ iff there exists 
an $S\subseteq A^n$ such that the ${\cal L}\cup\{X\}$--model 
${\goth A}_S=({\goth A},S)$, where $S$ is the interpretation of $X$ 
under ${\goth A}_S$, is a model of $\Gamma (X)$.

Let $^{*}V$ be a nonstandard universe. 
Let $\goth A$ be an $\cal L$--model with base set $A$ in the standard
universe.
We write $^{*}{\goth A}$ for an internally presented $\cal L$--model in 
$^{*}V$ with base set $^{*}A$ and the interpretation $P^{^{*}{\goth A}}=
\;^{*}(P^{\goth A})$ for every symbol $P\in {\cal L}$. 
It is not hard to see that
${\goth A}\equiv\; ^{*}{\goth A}$. In fact, $\goth A$ can be considered as
an elementary submodel of $^{*}{\goth A}$.

\bigskip

\noindent {\bf Main Theorem}\quad {\em Let $\kappa <\beth_{\omega}$ be
a regular cardinal. Then the following are equivalent:

(1) $I\!P_{\kappa}$,

(2) For any first--order language $\cal L$ with fewer than $\kappa$ 
many symbols, for any $n$--$\Delta^1_0 ({\cal L})$ type $\Gamma (X)$
and for any internally presented $\cal L$--model
$\goth A$ in $^{*}V$, if $\Gamma (X)$ is consistent with $\goth A$,
then $\goth A$ realizes $\Gamma (X)$.}

\bigskip

We will break the main theorem into following two theorems.

\begin{theorem}

Assume $\kappa<\beth_{\omega}$ is a regular cardinal. Let $^{*}V$ be a 
nonstandard universe which satisfies $I\!P_{\kappa}$. 
For any first--order language $\cal L$ with fewer than $\kappa$ many symbols,
for any $n$--$\Delta^1_0 ({\cal L})$ type $\Gamma (X)$
and for any internally presented $\cal L$--model
$\goth A$ in $^{*}V$, if $\Gamma (X)$ is consistent with $\goth A$,
then $\goth A$ realizes $\Gamma (X)$.

\end{theorem}

\noindent {\bf Proof:}\quad
Let $^{*}V$, $\cal L$, $\Gamma (X)$ and $\goth A$ are as described
in the theorem. We want to show that $\goth A$ realizes $\Gamma (X)$.

Since $\Gamma (X)$ is consistent with $\goth A$, there exists an
$\cal L$--model $\goth B$ with base set $B$ and an $S'\subseteq B^n$
such that the ${\cal L}\cup\{X\}$--model ${\goth B}_{S'}=({\goth B},S')$
is a model of $T\!h({\goth A})\cup\Gamma (X)$. We can assume $|B|\leq\kappa$ by
the Downward L\"{o}wenheim--Skolem--Tarski theorem. Furthermore we
can assume that ${\goth B}$ is in the standard universe because 
$\kappa<\beth_{\omega}$. Let $^{*}{\cal B}_{S'}=(^{*}{\cal B},\;^{*}S')$
be the internally presented ${\cal L}\cup\{X\}$--model in $^{*}V$ 
defined above.
It is easy to see now that ${\goth A}\equiv \;^{*}{\goth B}$.
By $I\!P_{\kappa}$, there is an isomorphism $i$ from $^{*}{\goth B}$
to $\goth A$. Let 
\[S=\{(i(b_1),i(b_2),\ldots,i(b_n)): (b_1,b_2,\ldots,b_n)\in\; ^{*}S'\}.\] 
Then $i$ is an isomorphism from $(^{*}{\goth B},\; ^{*}S')$ to
$({\goth A},S)$ in ${\cal L}\cup\{X\}$. 
Since ${\goth B}_{S'}\models\Gamma (X)$, then $^{*}{\goth B}_{S'}\models
\Gamma (X)$. Since $^{*}{\goth B}_{S'}\cong {\goth A}_{S}$, we conclude that
$\goth A$ realizes $\Gamma (X)$. \quad $\square$

\bigskip

\noindent {\bf Remark:}\quad 
If we replace the predicate symbol $X$ in the definition of $n$--$\Delta_0^1
({\cal L})$ types by a new constant symbol $c$, the proof of Theorem 1 can still
go through. So, as a corollary of Theorem 1, $I\!P_{\kappa}$ implies 
$\kappa$--saturation.

\bigskip

Next we are going to prove the converse of Theorem 1.
Before going further, we need to introduce more notation. The first is
about model pairs. Let $\cal L$ be a language. We call a language 
${\cal L}'$ {\em the language for $\cal L$--model pairs} if

(1) $\cal L$ and ${\cal L}'$ have same function symbols,  

(2) every relation symbol in $\cal L$ is in ${\cal L}'$ and ${\cal L}'$ 
contains two additional unary relation symbols $P$ and $Q$,

(3) for every constant symbol $c$ in $\cal L$, ${\cal L}'$ contains 
exactly two copies of $c$, say, $c_0$ and $c_1$.

Let $\goth A$ and $\goth B$ be two $\cal L$--models. A model pair
${\goth C}_{{\goth A},{\goth B}}$ is an ${\cal L}'$--model with
base set $A\cup B$ (we assume $A\cap B=\emptyset$) such that

(1) for every function symbol or relation symbol $R$ in $\cal L$,
$R^{{\goth C}_{{\goth A},{\goth B}}}=R^{\goth A}\cup R^{\goth B}$,

(2) $P^{{\goth C}_{{\goth A},{\goth B}}}=A$ and $Q^{{\goth C}_{{\goth A},
{\goth B}}}=B$,

(3) $c_{0}^{{\goth C}_{{\goth A},{\goth B}}}=c^{\goth A}$ and
$c_{1}^{{\goth C}_{{\goth A},{\goth B}}}=c^{\goth B}$.

Let $\cal L$ be a language and let $R$ be an unary predicate symbol.
For any $\cal L$--formula $\phi$, we write $\phi^{R}$, the relativization
of $\phi$ under $R$, for the formula defined inductively by 

(1) $\phi^{R}=\phi$ if $\phi$ is an atomic formula,

(2) if $\phi =\neg \psi$, then $\phi^{R}=\neg\psi^{R}$,

(3) if $\phi =\psi\wedge\chi$, then $\phi^{R}=\psi^{R}\wedge\chi^{R}$,

(4) if $\phi =\exists x\psi$, then $\phi^{R}=\exists x(R(x)\wedge\psi^{R})$.

\begin{theorem}

Let $^{*}V$ be a nonstandard universe. Let $\kappa$ be a regular cardinal.
If for any language $\cal L$ with fewer than $\kappa$--many symbols, 
for any internally presented $\cal L$--model $\goth A$ in $^{*}V$, 
and for any $2$--$\Delta^1_0 ({\cal L})$ type $\Gamma (X)$ which is consistent
with $\goth A$, the model $\goth A$ realizes $\Gamma (X)$, 
then $^{*}V$ satisfies $I\!P_{\kappa}$.

\end{theorem}

\noindent {\bf Proof:}\quad
Let $\cal L$ be a language with fewer than $\kappa$--many symbols.
Let $\goth A$ and $\goth B$ be two in ternally presented $\cal L$--models
in $^{*}V$ such that ${\goth A}\equiv {\goth B}$. We want to show that
${\goth A}\cong {\goth B}$.
    
Let ${\cal L}'$ be the language for $\cal L$--model pairs and let
${\goth C}_{{\goth A},{\goth B}}$ be the model pair of $\goth A$ and
$\goth B$. We want now to define a $2$--$\Delta^1_0 ({\cal L}')$ 
type $\Gamma (X)$ which will be used to force an isomorphism between $\goth A$
and $\goth B$. Let
\[\Gamma (X)=\{\phi_{n}(X):n=0,1,2,3,4\}\cup\{\psi_{\varphi}(X):
\varphi\mbox{ is an }\cal L\mbox{--formula.}\},\] where
\[\phi_0 (X)=\forall x\forall y (X(x,y)\rightarrow P(x)\wedge Q(y))\]
\[\phi_1 (X)=\forall x (P(x)\rightarrow\exists y X(x,y))\]
\[\phi_2 (X)=\forall y (Q(y)\rightarrow \exists x X(x,y))\]
\[\phi_3 (X)=\forall x\forall y\forall z (X(x,z)\wedge X(y,z)\rightarrow x=y)\]
\[\phi_4 (X)=\forall x\forall y\forall z (X(z,x)\wedge X(z,y)\rightarrow x=y)\]
\[\psi_{\varphi} (X)=\forall x_1\cdots\forall x_n\forall y_1\cdots\forall y_n
(\bigwedge_{k=1}^{n} X(x_k,y_k)\rightarrow (\varphi^{P}(x_1,\ldots,x_n)
\leftrightarrow\varphi^{Q}(y_1,\ldots,y_n))).\]
We can see that the sentences $\{\phi_n (X):n=0,1,2,3,4\}$ say that $X$
is a one to one correspondence between $P$ and $Q$. Hence $\Gamma (X)$
says that the one to one correspondence $X$ is actually an isomorphism
between $\goth A$ and $\goth B$.
It is easy to check that for any two $\cal L$--models ${\goth A}'$ and
${\goth B}'$, the model pair ${\goth C}_{{\goth A}',{\goth B}'}$ realizes
$\Gamma (X)$ if and only if ${\goth A}'\cong {\goth B}'$. We need now only
to show that $\Gamma (X)$ is consistent with ${\goth C}_{{\goth A},{\goth B}}$.
Since $\goth A$ and $\goth B$ are elementarily equivalent, there exists
an ultrafilter $\cal F$ on some cardinal $\lambda$ such that the ultrapower
of $\goth A$ and the ultrapower of $\goth B$ modulo $\cal F$ are isomorphic 
(see [S]). 
Hence the ultrapower of ${\goth C}_{{\goth A},{\goth B}}$ modulo $\cal F$,
which is the model pair of the ultrapower of $\goth A$ and the ultrapower
of $\goth B$ modulo $\cal F$,
realizes $\Gamma (X)$. On the other hand, the ultrapower of ${\goth C}_
{{\goth A},{\goth B}}$ is elementarily equivalent to ${\goth C}_
{{\goth A},{\goth B}}$. So $\Gamma (X)$ is consistent with ${\goth C}_
{{\goth A},{\goth B}}$. \quad $\square$

\bigskip

\noindent {\bf Remarks:}\quad (1) As a corollary we have that in a nonstandard
universe, the realizability for all $2$--$\Delta^1_0 (\cal L)$ types is
equivalent to the realizability of all $n$--$\Delta^1_0 (\cal L)$ types
for every $n$. (2) We didn't required that $\kappa<\beth_{\omega}$ in
Theorem 2.

\section{The applications}

The first application will give an answer to Ross's question about ($\dagger$).
In order to avoid dealing with the lengthy definition of Loeb measure we
are going to express ($\dagger$) in an internal version as Ross did (see [R]).

We use the words {\em finite} or {\em infinite} for externally finite or
externally infinite, respectively. We use $^{*}$finite or $^{*}$infinite
for internally finite or internally infinite, respectively. For example,
if $n\in\;^{*}{\Bbb N}\smallsetminus {\Bbb N}$, where ${\Bbb N}$ is
the set of all standard natural numbers, then the set $\{0,1,\ldots,n\}$
is both $^{*}$finite and infinite. We use $\Bbb R$ for the set of all 
standard reals.

Let $^{*}V$ be a nonstandard universe. Let $r\in\;^{*}{\Bbb R}$. We say that
$r$ is finite if there is a standard $n\in {\Bbb N}$ such that
$|r|<n$. Otherwise we call $r$ infinite. We say that $r$ is an infinitesimal 
if $|r|<\frac{1}{n}$ for every standard $n\in {\Bbb N}$. 

\begin{application}

($I\!P_{\aleph_0}$) Suppose $\Omega$ is an infinite internal set and
$\cal B$ is an internal subalgebra of $^{*}{\cal P}(\Omega)$ which contains
all singletons. Let $\mu :{\cal B}\rightarrow ^{*}[0,\infty)$ be an internal,
finitely additive measure with $\mu (\Omega)$ infinite and $\mu (\{x\})$
infinitesimal for every $x\in\Omega$. Then there exists a subset $S\subseteq
\Omega$ such that

(1) for any $D\in {\cal B}$ with $\mu (D)$ finite, for any $n\in {\Bbb N}$,
there exists an $E\in {\cal B}$ such that $D\cap S\subseteq E$ and
$\mu (E) <\frac{1}{n}$,

(2) for any $D\in {\cal B}$, if $S\subseteq D$, then $\mu (D)$ is infinite.

\end{application}

\noindent {\bf Proof:}\quad
Let $\Omega$, $\cal B$ and $\mu$ be as described in the Application 1.
Let $^{*}{\Bbb R}$ be the set of hyperreal numbers. Assume that
$\Omega$, $\cal B$ and $^{*}{\Bbb R}$ are all disjoint. 

We form first an 
internally presented ${\cal L}_{\goth A}$--model $\goth A$ with base
set $A=\Omega\cup {\cal B}\cup ^{*}{\Bbb R}$ such that
\[{\goth A}=(A\; ;\Omega,{\cal B}, ^{*}{\Bbb R},\in,\mu,\cap,\smallsetminus,
+,*,\leq,0,1),\]
where $\Omega$, $\cal B$ and $^{*}{\Bbb R}$ are three unary relations on $A$,
$\in\;\subseteq\Omega\times {\cal B}$ is the membership relation,
$\mu:{\cal B}\mapsto\;^{*}{\Bbb R}$ is the finite additive measure,
$\cap$ is the set intersection and $\smallsetminus$ is the set subtraction
on $\cal B$, and $( ^{*}{\Bbb R}\; ;+,*,\leq,0,1)$ is the usual hyperreal
ordered field. For simplicity, we do not distinguish a symbol in 
${\cal L}_{\goth A}$ from its interpretation under $\goth A$.

We form next a $1$--$\Delta^1_0 ({\cal L}_{\goth A})$ type 
$\Gamma (X)$ such that
\[\Gamma (X)=\{\phi (X),\psi_n (X),\chi_n (X): n=1,2,\ldots\},\]
where \[\phi (X)=\forall x(X(x)\rightarrow\Omega(x))\]
\[\psi_n (X)=\forall U({\cal B}(U)\wedge\mu (U)<n\rightarrow
\exists V({\cal B}(V)\wedge\forall x(X(x)\wedge x\in U\rightarrow 
x\in V)\wedge\mu (V) <\frac{1}{n}))\]
\[\chi_n (X)=\forall U({\cal B}(U)\wedge\forall x(X(x)\rightarrow x\in U)
\rightarrow \mu (U)>n).\]
Notice that in $\goth A$, the element $1$ is definable, so do $n$ and $\frac{1}{n}$
for every $n\in {\Bbb N}$.

The sentence $\phi (X)$ says that $X$ is a subset of $\Omega$. The sentence
$\psi_n (X)$ says that the intersection
of $X$ with any $U$ in $\cal B$ with measure less than $n$ has outer measure
less than $\frac{1}{n}$. The sentence $\chi_n (X)$ says that $X$ has outer
measure greater than $n$.
So the application 1 is true if and only if $\goth A$ realizes $\Gamma (X)$.
Hence, by the Theorem 1, it suffices to show that $\goth A$ is consistent
with $\Gamma (X)$.

Let $T=T\!h(\goth A)$.

\bigskip

{\bf Claim:}\quad $T\cup\Gamma (X)$ is consistent.

Proof of Claim:\quad 
By Downward L\"{o}wenheim-Skolem Theorem we can find a countable
model ${\goth A}_0\preccurlyeq {\goth A}$ 
with base set $A_0=\Omega_0\cup {\cal B}_0
\cup {\Bbb R}_0$. Since \[\exists U({\cal B}(U)\wedge\forall x(\Omega(x)
\rightarrow x\in U))\] is true in $\goth A$, it is true in ${\goth A}_0$.
Hence $\Omega_0\in {\cal B}_0$. Since $\mu (\Omega)>n$ for all
$n\in {\Bbb N}$ are true in $\goth A$, they are also true in ${\goth A}_0$.
Hence $\mu (\Omega_0)$ is infinite in ${\goth A}_0$.
Since \[\forall U\forall x\forall y({\cal B}(U)\wedge\Omega(x)\wedge\Omega(y)
\wedge(x\in U\wedge y\in U\rightarrow x=y)\rightarrow\mu (U)<\frac{1}{n}\]
for all $n\in {\Bbb N}$ are true in $\goth A$, they are also true in 
${\goth A}_0$.
Hence the measure of every singleton is infinitesimal in ${\goth A}_0$.
Let \[\{B\in{\cal B}_0:\mu (B)\mbox{ is finite }\}=\{B_n:n\in {\Bbb N}\}.\]
It is now easy to pick \[x_n\in\Omega_0\smallsetminus(\bigcup_{k=0}^{n-1}B_k
\cup\{x_k:k<n\})\] because $\Omega_0$ has infinite measure and the measure
of $\bigcup_{k=0}^{n-1}B_k\cup\{x_k:k<n\}$ is finite. Also notice that
the measure of a finite set $\{x_k:k<n\}$ for $n\in {\Bbb N}$ is 
infinitesimal because the sum of finitely many infinitesimals is 
an infinitesimal and ${\cal B}_0$ is closed under finite union.

Let $S_0=\{x_n:n\in {\Bbb N}\}$. It is obvious that $({\goth A}_0,S_0)$
is a model of $T\cup\Gamma (X)$. \quad $\square$

\bigskip

Next three applications are also the questions of [R] and were proved in [J].
The purpose of including them here with simplified proofs is to illustrate 
that $I\!P_{\kappa}$ is an ``easy to use'' tool in nonstandard analysis. 

\begin{application}

($I\!P_{\aleph_0}$) Suppose that $(P,<_{P})$ and $(Q,<_{Q})$ are two
internal linear orders without endpoints. There is an order--preserving
map $f:P\mapsto Q$ such that $f[P]$ is cofinal in $Q$.

\end{application}

\noindent {\bf Proof:}\quad
Without loss of generality, we can assume that $P\cap Q=\emptyset$.
Let $\goth A$ be an internally presented ${\cal L}_{\goth A}$--model with
base set $A=P\cup Q$
such that \[{\goth A}=(A\; ;P,Q,\leq_{P},\leq_{Q}),\]
where $P$ and $Q$ are two binary relations on $A$, and $<_{P}$ and $<_{Q}$
are the correspondent orders on $P$ and $Q$. We define a $2$--$\Delta_0^1 
({\cal L}_{\goth A})$ type $\Gamma (X)$ such that
\[\Gamma (X)=\{\phi (X),\psi (X),\chi (X),\pi (X)\}\]
where \[\phi (X)=\forall x\forall y (X(x,y)\rightarrow P(x)\wedge Q(y))\]
\[\psi (X)=\forall x\exists ! y(P(x)\rightarrow X(x,y))\]
\[\chi (X)=\forall x_0\forall x_1\forall y_0\forall y_1
(x_0 <_{P}x_1\wedge X(x_0,y_0)\wedge X(x_1,y_1)\rightarrow y_0 <_{Q} y_1)\]
\[\pi (X)=\forall y_0\exists x\exists y_1 (Q(y_0)\rightarrow P(x)\wedge
X(x,y_1)\wedge y_0 <_{Q} y_1).\]
In $\Gamma (X)$ the sentence $\phi (X)$ says that $X$ is a relation
between $P$ and $Q$, the sentence $\psi (X)$ says that $X$ is the graph 
of a function from $P$ to $Q$, the sentence $\chi (X)$ says that
the function is order--preserving, and $\pi (X)$ says that the function
is a cofinal embedding.
If there exists an $S\subseteq P\times Q$ such that $({\goth A},S)\models
\Gamma (X)$, then it is easy to see that the map $f$ defined by its graph
$S$ is the order--preserving map we are looking for. By Theorem 1, we need
only to show that $T\cup\Gamma (X)$ is consistent, where $T=T\!h(\goth A)$.

\bigskip

{\bf Claim:}\quad $T\cup\Gamma (X)$ is consistent.

Proof of Claim:\quad
Let ${\goth A}_0=(A_0;P_0,Q_0,\leq_{P_0}\leq_{Q_0})$ be a countable 
elementary submodel of $\goth A$.
By the Compactness Theorem and L\"{o}wenheim--Skolem Theorem 
${\goth A}_0$ can be elementarily extended to a countable model
${\goth A}_1=(A_1,P_1,Q_1,\leq_{P_1},\leq_{Q_1})$ such that the set of
all rational numbers in $[0,1)$, together with the usual order, can be
order-isomorphically embedded into the set $\{q\in Q_1:\forall x\in Q_0\;
(x\leq_{Q_1} q)\}$. Let $i_0$ be that embedding. By the same argument
we can find an elementary chain of length $\omega$ of countable models
${\goth A}_0\preccurlyeq {\goth A}_1\preccurlyeq\cdots$ and a sequence of maps
$\{i_n:n\in\omega\}$ such that $i_n$ is an order--preserving map
from the set of all rational numbers in $[n, n+1)$ with the usual
order to $\{q\in Q_{n+1}:\forall x\in Q_n\;(x\leq_{Q_{n+1}}q)\}$.
Let ${\goth A}_{\omega}=\bigcup_{n\in\omega}{\goth A}_n$ and let
$i_{\omega}=\bigcup_{n\in\omega}i_n$. Since ${\goth A}_0$ is elementarily
equivalent to both $\goth A$ and ${\goth A}_{\omega}$, then 
${\goth A}_{\omega}$ is a model of $T$.
It is easy to see that $i_{\omega}$ is an
order-preserving map from the set of all positive rational numbers
cofinally into $Q_{\omega}$. Since every countable order without a right
endpoint can be cofinally embedded into the set of all 
positive rational numbers,
then $P_{\omega}$ can be cofinally embedded into $Q_{\omega}$. Hence
${\goth A}_{\omega}$ realizes $\Gamma (X)$. This proves the consistency
of $T\cup\Gamma (X)$.\quad $\square$

\begin{application}

($I\!P_{\aleph_0}$) Let $(P,<_{P})$ be an internal partial order with
no right endpoints. There is an external subset $S\subseteq P$ such that
$\{s\in S:s<_{P}p\}$ is internal for every $p\in P$.

\end{application}

\noindent {\bf Proof:}\quad Let ${\cal Q}=\; ^{*}{\cal P}(P)$.
Let $\goth A$ be an internally presented ${\cal L}_{\goth A}$--model
with base set $A=P\cup {\cal Q}$ such that
\[{\goth A}=(A\; ; P,{\cal Q},\in,<_{P},\cap,\smallsetminus),\]
where $P$ and $\cal Q$ are two unary relations, $\in \;\subseteq 
P\times {\cal Q}$ is the membership relation, $<_{P}$ is the order on $P$,
$\cap$ is the set intersection on $\cal Q$ and $\smallsetminus$ is
the set subtraction on $\cal Q$.
We define a $1$--$\Delta_0^1 ({\cal L}_{\goth A})$ type $\Gamma (X)$ such that
\[\Gamma (X)=\{\phi (X),\psi (X),\chi (X)\}\] where
\[\phi (X)=\forall x(X(x)\rightarrow P(x))\]
\[\psi (X)=\forall x\exists U (P(x)\rightarrow {\cal Q}(U)\wedge\forall y
(y<_{P}x\wedge X(y)\leftrightarrow y\in U))\]
\[\chi (X)=\forall U\exists x({\cal Q}(U)\rightarrow 
(x\in U\wedge\neg X(x))\vee (X(x)\wedge\neg x\in U)).\]

The sentence $\phi (X)$ says that $X$ is a subset of $P$, the sentence 
$\psi (X)$ says that for every $x$ in $P$ there exists a $U$ in 
$^{*}{\cal P}(P)$ such that $U=\{y\in X:y<_{P}x\}$, and the sentence
$\chi (X)$ says that for all $U$ in $^{*}{\cal P}(P)$, the set $X$
is different from $U$.
It is easy to see that if there is an $S\subseteq P$ such that
$({\goth A},S)\models\Gamma (X)$, then $S$ the set we are looking for.
By Theorem 1, it suffices to show that $T\cup\Gamma (X)$ is consistent,
where $T=T\!h(\goth A)$.

\bigskip

{\bf Claim:}\quad $T\cup\Gamma (X)$ is consistent.

Proof of Claim:\quad
Let ${\goth A}_0=(A_0;P_0,{\cal Q}_0,\ldots)$ be a countable elementary 
submodel of $\goth A$. It suffices to construct a set $S=\{s_n:n\in\omega\}
\subseteq P_0$ such that $({\goth A}_0,S)\models\Gamma (X)$.
Let $P_0=\{p_n:n\in\omega\}$ and let ${\cal Q}_0=\{Q_n:n\in\omega\}$.
Since $P$ has no right endpoints, then $P_0$ has no right endpoints.
Now we can pick the elements $s_k$ and $t_k$ from $P_0$
for every $k\in\omega$ such that 

(1) $s_0<t_0<s_1<t_1<\cdots$,

\noindent and for every $k\in\omega$ 

(2) we have $s_k\not<p_k$ and 

(3) either both
$s_k$ and $t_k$ are in $Q_k$ or both $s_k$ and $t_k$ are not in $Q_k$.

\noindent Let $S=\{s_n:n\in\omega\}$. Then $S$ differs from every element
in ${\cal Q}_0$. For every $p\in P_0$ the set $\{s\in S:s\leq p\}$
is finite, and hence is in ${\cal Q}_0$ because for every finite
set $\{a_1,\ldots,a_n\}\subseteq P_0$ the sentence
\[\exists U({\cal Q}(U)\wedge\forall x(x\in U\leftrightarrow
\bigvee_{i=1}^{n}x=a_i))\]
is true in $\goth A$ and therefore, it is true in ${\goth A}_0$.

The arguments above showed that $({\goth A}_0,S)\models\Gamma (X)$.
\quad $\square$

\begin{application}

($I\!P_{\aleph_0}$) Let $P$ and $Q$ be two $^{*}$infinite internal sets
There is a bijection $f:P\mapsto Q$ such that for every $^{*}$finite
$b\subseteq P$ and for every $^{*}$finite $c\subseteq Q$, the restriction
of $f$ to $b$ and the restriction of $f^{-1}$ to $c$ are internal.

\end{application}

\noindent {\bf Proof:}\quad
Without loss of generality, we can assume that $P\cap Q=\emptyset$.
Let $\goth A$ be an internally presented ${\cal L}_{\goth A}$--model
with base set $A=P\cup Q\cup F$, where $F=\{f:f$ is an internal
bijection from some $^{*}$finite subset of $P$ to $Q\}$, such that
\[{\goth A}=(A\; ;P,Q,F,R)\]
where $P$, $Q$ and $F$ are three unary relations on $A$ and
$R\subseteq P\times Q\times F$ is defined by
\[(a,b,f)\in R\;\;\mbox{ iff }\;\;(a,b)\in f.\]

We now define a $2$--$\Delta_0^1 ({\cal L}_{\goth A})$ type $\Gamma (X)$ such 
that
\[\Gamma (X)=\{\phi_n (X):n=0,1,2,3,4,5\},\]
where \[\phi_0 (X)=\forall x\forall y (X(x,y)\rightarrow P(x)\wedge Q(y))\]
\[\phi_1 (X)=\forall x\forall y\forall z (X(x,z)\wedge X(y,z)\rightarrow x=y)\]
\[\phi_2 (X)=\forall x\forall y\forall z (X(z,x)\wedge X(z,y)\rightarrow x=y)\]
\[\phi_3 (X)=\forall x\exists y ((P(x)\rightarrow X(x,y))\wedge
(Q(x)\rightarrow X(y,x)))\]
\[\phi_4 (X)=\]\[\forall g\exists f(F(g)\rightarrow F(f)\wedge
\forall x(\exists y R(x,y,g)\leftrightarrow\exists y R(x,y,f))\wedge
\forall x\forall y (R(x,y,f)\rightarrow X(x,y)))\]
\[\phi_5 (X)=\]\[\forall g\exists f(F(g)\rightarrow F(f)\wedge
\forall y(\exists x R(x,y,g)
\leftrightarrow\exists x R(x,y,f))\wedge\forall x\forall y (R(x,y,f)
\rightarrow X(x,y))).\]
The sentences $\phi_0 (X),\phi_1 (X),\phi_2 (X),\phi_3 (X)$ 
say that $X$ is a one to one
onto correspondence between $P$ and $Q$. The sentence $\phi_4 (X)$
says that the restriction of $X$ on any $^{*}$finite set of $P$ (as the
domain of an element $g$ in $F$) coincides with
an element $f$ in $F$. The sentence $\phi_5 (X)$ says that the restriction of
$X$ on any $^{*}$finite subset of $Q$ (as the range of an element $g$ in $F$) 
coincides also with an element $f$ in $F$.
It is easy to see that if there exists an $S\subseteq P\times Q$
such that $({\goth A},S)\models\Gamma (X)$, then the bijection induced by
$S$ is the map we are looking for. Let $T=T\!h({\goth A})$.
By Theorem 1, we need only to show that 

\bigskip

{\bf Claim:}\quad $T\cup\Gamma (X)$ is consistent.

Proof of Claim:\quad
Let ${\goth A}_0=(A_0;P_0,Q_0,F_0,R_0)$ be a countable elementary
submodel of $\goth A$. It suffices to find a relation
$S\subseteq P_0\times Q_0$ such that $S$ is the graph of a bijection
$i$ from $P_0$ to $Q_0$ and for every $C$ which is the domain of a function in
$F_0$ and for every $D$ which is the range of a function in $F_0$, both
$i\!\res\!C$ and $(i^{-1}\!\res\!D)^{-1}$ are functions in $F_0$.

Let $F_0=\{f_n:n\in\omega\}$. We want to construct a sequence 
$\{i_m:m\in\omega\}\subseteq F_0$ such that 

(1) $i_0\subseteq i_1\subseteq \i_2\subseteq\cdots$, 

(2) for every $f\in F_0$ there is an $m\in\omega$
such that $dom(f)\subseteq dom(i_m)$ and $range(f)\subseteq range(i_m)$.

The claim follows from the construction because we can let 
$i=\bigcup_{m\in\omega}i_m$. It is easy to check that

(a) $i$ is one to one function, 

(b) $dom(i)=P_0$ and $range(i)=Q_0$,

(c) for every $f\in F_0$ there exists an $m\in\omega$ such that
$i\!\res\!dom(f)=i_m\!\res\!dom(f)\in F_0$ and $(i^{-1}\!\res\!range(f))^{-1}
=(i_m^{-1}\!\res\!range(f))^{-1}\in F_0$.

\noindent (c) is true because both sentences
\[\forall g\in F\forall f\in F(dom(g)\subseteq dom(f)\rightarrow\exists
h\in F(h=f\!\res\!dom(g)))\] and
\[\forall g\in F\forall f\in F(range(g)\subseteq range(f)\rightarrow\exists
h\in F(h=(f^{-1}\!\res\!range(g))^{-1}))\] are true in $\goth A$.
Then let $S$ be the graph of $i$ and we have now $({\goth A}_0,S)\models
\Gamma (X)$.

Before constructing $\{i_m:m\in\omega\}$ let's observe two facts.

{\em Fact 1}. Suppose $f, g\in F_0$. If $dom(f)\cap dom(g)=\emptyset$ and
$range(f)\cap range(g)=\emptyset$, then $f\cup g\in F_0$.

{\em Fact 2}. For any $f,g\in F_0$ there are $h,j\in F_0$ such that 
$dom(f)=dom(h)$, $range(g)\cap range(h)=\emptyset$ and $range(f)
=range(j)$, $dom(g)\cap dom(j)=\emptyset$.

\noindent {\em Fact 1} and {\em Fact 2} above are true in 
${\goth A}_0$ because they are true in $\goth A$.

Let $i_0=f_0$. Assume that we have constructed $\{i_m:m<k\}\subseteq
F_0$ such that \[i_0\subseteq i_1\subseteq\cdots\subseteq i_{k-1},\]
\[dom(f_m)\subseteq dom(i_{2m}) \mbox{ when }2m<k\] 
and \[range(f_m)\subseteq
range(i_{2m+1})\mbox{ when }2m+1<k.\]

Case 1:\quad $k=2n$.

Let $C=dom(f_n)\smallsetminus dom(i_{k-1})$. If $C=\emptyset$, then
let $i_{k}=i_{k-1}$. Otherwise let $h\in F_0$ such that $dom(h)=C$
and $range(h)\cap range(i_{k-1})=\emptyset$. The function $h$ exists
by {\em Fact 2}. Let $i_{k}=i_{k-1}\cup h$. The function $i_k\in F_0$
by {\em Fact 1}.

Case 2:\quad $k=2n+1$.

Let $D=range(f_n)\smallsetminus range(i_{k-1})$. If $D=\emptyset$,
then let $i_{k}=i_{k-1}$. Otherwise let $h\in F_0$ such that
$range(h)=D$ and $dom(h)\cap dom(i_{k-1})=\emptyset$. Again the function
$h$ exists by {\em Fact 2} and $i_{k}=i_{k-1}\cup h\in F_0$ by 
{\em Fact 1}.
\quad $\square$

\bigskip

\noindent {\bf Remark}:\quad The claims in above four applications
could also be shown by quoting simply four results from [R] or other
papers.
We present our own proofs here because these proofs use only countable
models so that the reader can read the paper without knowing special
models.

\bigskip

Department of Mathematics,
University of California, Berkeley, CA 94720, USA.

\bigskip

Institute of Mathematics, 
The Hebrew University, 
Jerusalem, Israel.

\bigskip

Department of Mathematics, 
Rutgers University, 
New Brunswick, NJ, 08903, USA.

\bigskip

{\em Sorting:} The first address is the first author's; the last
two are the second author's.

\end{document}